\documentclass[11pt]{article}
\usepackage{amsmath,amssymb,amsthm}
\usepackage{mathrsfs}
\usepackage{enumerate}
\usepackage{enumitem}
\usepackage{latexsym,amsmath,amssymb,amsfonts,epsfig,graphicx,cite,psfrag}
\usepackage{eepic,color,colordvi,amscd}
\usepackage{ebezier}
\usepackage{verbatim}

\newtheorem{thm}{Theorem}[section]
\newtheorem{lm}[thm]{Lemma}

\newtheorem{conj}{Conjecture}[section]
\newtheorem{pro}[thm]{Proposition}

\usepackage{subfigure}

\newcommand{\ar}{{\rm ar}}
\newcommand{\ex}{{\rm ex}}
\newcommand{\G}{{\cal G}}

\begin{document}
	
\title{Anti-Ramsey number of matchings in $3$-uniform hypergraphs}

\author{Mingyang Guo\footnote{School of Mathematics and Statistics,
Xi'an Jiaotong University, Xi'an, Shaanxi 710049, China. Email: {\tt g1556010753g@163.com}.}
\and
Hongliang Lu \footnote{School of Mathematics and Statistics,
Xi'an Jiaotong University, Xi'an, Shaanxi 710049, China. Email: {\tt luhongliang@mail.xjtu.edu.cn},
Research is supported in part by  National Natural Science Foundation of China  (No.\ 12271425).}
\and Xing Peng\thanks{Center for Pure Mathematics, School of Mathematical Sciences, Anhui University,  Hefei, Anhui 230601, China. E-mail: {\tt x2peng@ahu.edu.cn}, Research is supported in part by National Natural Science Foundation of China (No.\ 12071002) and the Anhui Provincial Natural Science Foundation (No. 2208085J22).}
}

%
%
%

\date{}

\maketitle

\date{}

\maketitle

\begin{abstract}
Let $n,s,$ and $k$ be positive integers such that $k\geq 3$, $s\geq 3$ and $n\geq ks$. An $s$-matching $M_s$ in a $k$-uniform hypergraph is a set of $s$ pairwise disjoint edges. The anti-Ramsey number $\ar(n,k,M_s)$ of an $s$-matching  is the smallest integer $c$ such that each edge-coloring of the $n$-vertex $k$-uniform complete hypergraph with exactly $c$ colors contains an $s$-matching  with distinct colors. In 2013, \"Ozkahya and Young proposed a conjecture on the
 exact value of $\ar(n,k,M_s)$  for all  $n \geq sk$ and $k \geq 3$. A 2019 result by Frankl and Kupavskii  verified this conjecture for all $n \geq sk+(s-1)(k-1)$ and $k \geq 3$.
  We aim to determine the value of $\ar(n,3,M_s)$ for $3s \leq n < 5s-2$ in this paper. Namely,  we prove that
 if $3s<n<5s-2$ and $n$ is large enough, then $\ar(n,3,M_s)=\ex(n,3,M_{s-1})+2$.
Here  $\ex(n,3,M_{s-1})$ is the Tur\'an number of an $(s-1)$-matching.
 Thus this result confirms the conjecture of \"Ozkahya and Young for $k=3$, $3s<n<5s-2$ and sufficiently large $n$.
   For $n=ks$ and $k\geq 3$, we present a new construction for the lower bound of $\ar(n,k,M_{s})$ which shows the conjecture by \"Ozkahya and Young is not true. In particular, for $n=3s$, we prove that $\ar(n,3,M_s)=\ex(n,3,M_{s-1})+5$ for sufficiently large $n$.
  \end{abstract}

\section{Introduction}
Given a graph $G$, the {\it anti-Ramsey number} $\ar(n,G)$  is the smallest integer $c$ such that each edge-coloring of $K_n$ with exactly $c$ colors contains a  rainbow copy of $G$. Here $G$ is rainbow if  all edges have distinct colors. For a set of graphs $\G$, the Tur\'an number $\ex(n,\G)$ is the maximum possible number of edges in an $n$-vertex graph which does not contain any $H \in \G$ as a subgraph.
The value of $\ar(n,G)$ is closely related to the Tur\'an number $\ex(n,G)$ as the following inequality, see \cite{ESS},
\begin{equation} \label{lb}
2+\ex(n,\G) \leq \ar(n,G) \leq 1+\ex(n,G),
\end{equation}
where $\G=\{G-e: e\in E(G)\}$. For the lower bound, let $F$ be a $\G$-free graph with $\ex(n,\G)$ edges. One can assign distinct colors to edges of $F$ and one more color to edges of $\overline{F}$. It is easy to see that there is not rainbow $G$.
A seminal result by  Erd\H{o}s-Simonovits-S\'os \cite{ESS} asserts that $\ar(n,K_p)=\ex(n,K_{p-1})+2$ for $n$ large enough.  Montellano-Ballesteros \cite{MN} and Neumann-Lara \cite{OY} extended this result to all values of $n$ and $p$ with $n>p\geq3$. An $s$-{\it matching} $M_s$ is a set of $s$ independent edges. The Tur\'an number $\ex(n,M_s)$ was determined by Erd\H os and Gallai \cite{EG} for $n\geq 2s$ and $s\geq 1$. For the anti-Ramsey number of $M_s$, Schiermeyer  \cite{S}  first proved that $\ar(n,M_s)=\ex(n,M_{s-1})+2$ for $s\geq 2$ and $n\geq 3s+3$. Later, Fujita, Kaneko, Schiermeyer and Suzuki \cite{FKSS} established the same result for  $s\geq 2$ and $n\geq 2s+1$. Finally,  Chen, Li and Tu \cite{CLT} determined the exact value of $\ar(n,M_s)$ for all $s\geq 2$ and $n\geq 2s$.
Haas and Young \cite{HY} reproved the case of $n=2s$   by a simple argument.  There is a large volume of literature on the anti-Ramsey number of graphs. Interested readers are referred to the survey by Fujita, Magnant, and Ozeki \cite{FMO}.

The anti-Ramsey number of  hypergraphs can be defined similarly.
A \emph{hypergraph} $H$ is a pair $H=(V(H),E(H))$, where $V(H)$ is a set of vertices and $E(H)$ is a set of non-empty subsets of $V$. A hypergraph is \emph{$k$-uniform} if $E(H) \subseteq\binom{V}{k}$, where $\binom{V}{k}=\{T\subseteq{V}: |T|=k\}$. A $k$-uniform hypergraph is also called a \emph{$k$-graph} for convenience. Throughout this paper, we often identify $E(H)$ with $H$ when there is no confusion. In an edge-coloring of a hypergraph $H$, a subgraph $G\subseteq H$ is \emph{rainbow} if all edges of $G$ have distinct colors.
For a $k$-uniform hypergraph $G$,  the anti-Ramsey number $\ar(n,k,G)$ is the smallest integer $c$ such that each edge-coloring of the $n$-vertex $k$-uniform complete hypergraph with exactly $c$ colors contains a rainbow copy of $G$. The {\it Tur\'an number} $\ex(n,k,G)$ of a $k$-uniform hypergraph $G$ is the maximum possible number of edges in an $n$-vertex $k$-uniform hypergraph which does not contain $G$ as a subgraph. One can easily see that inequality \eqref{lb} can be extended to hypergraphs as follows:
\begin{equation} \label{lbhyper}
2+\ex(n,k,\G) \leq \ar(n,k,G) \leq 1+\ex(n,k,G),
\end{equation}
where $\G=\{G-e: e\in E(G)\}$.
 An  $s$-\emph{matching} $M_s$ in a $k$-uniform hypergraph $H$ is a set of pairwise disjoint edges. The number of edges in a matching $M_s$, denoted by $|M_s|$, is called the \emph{size} of the matching.
The size of the largest matching in $H$ is denoted by $\nu(H)$, known as the\emph{ matching number} of $H$. A matching is \emph{perfect} if it covers all vertices of $V(H)$. A matching with $s$ edges is called an \emph{$s$-matching}. For the anti-Ramsey number of hypergraph matchings, \"Ozkahya and Young \cite{OY} proposed the following conjecture.
\begin{conj}\label{mainconj}
	Let $k\geq 3$ and  $s\geq 3$. If $n>ks$, then $\ar(n,k,M_s)=\ex(n,k,M_{s-1})+2$. In addition, if $n=ks$, then
	$$\ar(n,k,M_s)=\begin{cases}
		\ex(n,k,M_{s-1})+2, & \text{if $s<c_k$};\\
		\ex(n,k,M_{s-1})+k+1, &\text{if $s\geq c_k$,}
	\end{cases}$$
where $c_k$ is a constant depending on $k$.
\end{conj}
 The exact value of the Tur\'an number of hypergraph matchings is still unknown in general. In 1965, Erd\H os \cite{E65} asked for the determination of the maximum possible number of edges that can appear in any $k$-graph $H$ with $\nu(H)\leq s$.  We next introduce two constructions.
 Let $n,s,$ and $k$ be three positive integers such that $k\geq2$ and $n\geq ks+k-1$. For any $U\subseteq [n]$ with $|U|=k(s+1)-1$, define $D^k_{n,s}(U)$ as a hypergraph with edge set $\{e\in \binom{[n]}{k}:e\subseteq U\}$. Let $U,W$ be a partition of $[n]$ such that $|W|=s$. Define $H^k_{n,s}(U,W)$ as a hypergraph with edge set $\{e\in\binom{[n]}{k}: e\cap W\neq\emptyset\}$.
When there is no confusion, we denote  $H^k_{n,s}(U,W)$ and $D^k_{n,s}(U)$ by $H^k_{n,s}$ and $D^k_{n,s}$, respectively.
 Clearly, $|E(H^k_{n,s})|=\binom{n}{k}-\binom{n-s}{k}$ and $|E(D^k_{n,s})=\binom{k(s+1)-1}{k}$.
Furthermore,  the matching number of $H^k_{n,s}(U,W)$ and $D^k_{n,s}(U)$ is $s$. Based on these two constructions,  Erd\H{o}s made the following conjecture.

\begin{conj}[Erd\H os Matching Conjecture \cite{E65}]\label{conj}
	For $n\geq ks$, $k\geq 2$ and $s\geq 1$,
	\begin{equation*}
		\ex(n,k,M_{s+1})=\max\left\{\binom{n}{k}-\binom{n-s}{k},\binom{k(s+1)-1}{k}\right\}.
	\end{equation*}
\end{conj}
We next discuss the regimes in which the maximum is achieved by the construction $D^k_{n,s}(U)$. Let  $s_0(n, k)$ be the smallest $s$ for which $\binom{n}{k}-\binom{n-s}{k} \leq \binom{k(s+1)-1}{k}$.
It is easy to see that
$$
\lim _{n \rightarrow \infty} \frac{s_0(n, k)}{n}=\alpha_k
$$
 is the solution of the equation
$$
1-\left(1-\alpha_k\right)^k=k^k \alpha_k^k,
$$
where $\alpha_k \in(0,1)$.
One can check that for all $k \geqslant 3$, we have
$
\tfrac{1}{k}-\frac{1}{2 k^2}<\alpha_k<\tfrac{1}{k}-\frac{2}{5 k^2}.
$

 In 2019, Frankl and Kupavskii \cite{FK19} proved a stability result on Erd\H os Matching Conjecture for $k\geq3$ and either $n\geq(s+\max\{25,2s+2\})k$ or $n\geq(2+o(1))sk$. The Erd\H{o}s Matching Conjecture  was settled  for $k = 3$ and sufficiently large $n$ in \cite{LM}. The case of $k=3$ was   completely resolved  in \cite{F17}. For the state of art of Erd\H{o}s Matching conjecture, readers are referred to \cite{AFHRRS,BDE76,EG,E65,ESS,Fr13,F17,Fr172,FK18,FK19,FLM12,FRR,HLS12,LM,RR,Zhao}.

 Conjecture \ref{mainconj} was verified by Frankl and Kupavskii \cite{FK19} for all $n \geq sk+(s-1)(k-1)$ and $k \geq 3$.
\begin{thm}[Frankl and Kupavskii \cite{FK19}]\label{2sk}
	For $n\geq sk+(s-1)(k-1)$ and $k\geq 3$, we have $\ar(n,k,M_s)=\ex(n,k,M_{s-1})+2=\binom{n}{k}-\binom{n-s+2}{k}+2$.
\end{thm}
The case where $sk \leq n <sk+(s-1)(k-1)$ is unsolved and we study the case of $k=3$ in this paper.  We will prove the following theorem.
\begin{thm}\label{mainthm}
	For sufficiently large $n$, the following holds
	$$\ar(n,3,M_s)=\begin{cases}
		\ex(n,3,M_{s-1})+2, & \text{if $3s<n<5s-2$};\\
		\ex(n,3,M_{s-1})+5, &\text{if $n=3s$}.
	\end{cases}$$
\end{thm}
Thus Theorem \ref{mainthm} proves Conjecture \ref{mainconj} for $k=3$ and $3s<n<5s-2$ with $n$ large enough. Moreover, Theorem \ref{mainthm} shows that Conjecture \ref{mainconj} is not true for $n=3s$ and $n$ large enough.
 There are some related results on the anti-Ramsey number of hypergraph matchings, for example, those from Jin \cite{J} as well as Xue, Shan and Kang \cite{XSK}.

The rest of the paper is organized as follows. In Section 2, we will introduce notation and several previous results needed for our proofs. In Section 3, we will present a new construction for  the lower bound on $\ar(n,k,M_s)$.
  In section 4,  we will prove a stability result on matchings in $3$-graphs, which is a key ingredient in the proof of Theorem \ref{mainthm}. We will prove   Theorem \ref{mainthm} in Section 5. In Section 6, we will mention a few concluding remarks.

\section{Preliminaries}
We collect several previous results and introduce necessary definitions in this section.
We use $\ell$-set to denote a set of $\ell$ elements.
For an integer $n$, let $[n]=\{1,2,\ldots,n\}$. For a vertex $v\in V(H)$, let $N_H(v)=\{f\in\binom{V}{k-1}:f\cup \{v\}\in E(H)\}$ and  $d_H(v)=|N_H(v)|$.
Given a $k$-graph $H$ we write $e(H)=|E(H)|$. For a function $h$ defined over $E(H)$ and $P\subseteq E(H)$, let $h(P)=\{h(e)\ |\ e\in P\}$. We first recall the following result on Erd\H os Matching Conjecture.
\begin{thm}[ \cite{FLM12,F17}]\label{erdos3}
	For $s\geq 2$ and $n\geq 3s$,
		\begin{equation*}
		\ex(n,3,M_s)=\max\left\{\binom{n}{3}-\binom{n-s+1}{3},\binom{3s-1}{3}\right\}.
	\end{equation*}
\end{thm}

For a $3$-graph $H$, let $\delta_1(H)=\min\{d_H(x):x\in V(H)\}$. For the vertex degree version of Erd\H{o}s Matching Conjecture, H\`{a}n, Person and Schacht \cite{HPS09} showed that for a $3$-graph $H$, $\delta_1(H)>(5/9+o(1)){|V(H)|\choose 2}$
is sufficient for the existence of a perfect matching of $H$.
K\"{u}hn, Osthus and Treglown \cite{KOT} proved the following stronger result.
\begin{thm}[K\"uhn, Osthus and Treglown \cite{KOT}]\label{degree-version}
	There exists an integer $n_0\in\mathbb{N}$ such that if $H$ is a $3$-graph with $n\geq n_0$ vertices, $s$ is an integer with $1\leq s\leq n/3$, and
	\begin{equation*}
		\delta_1(H)>\binom{n-1}{2}-\binom{n-s}{2},
	\end{equation*}
	then $\nu(H)\geq s$.
\end{thm}

When $n\equiv0\pmod{3}$ and $n\geq n_0$, if $\delta_1(H)>\binom{n-1}{2}-\binom{2n/3}{2}$, then $H$ has a perfect matching. This result  was proved independently by Khan \cite{K}.

%

Given two $k$-graphs $H_1, H_2$ and a real number $\varepsilon>0$, we say that $H_2$ \textit{$\varepsilon$-contains}  $H_1$ if $V(H_1) = V(H_2)$ and $|E(H_1)\backslash E(H_2)|\leq\varepsilon|V(H_1)|^k$. In particular,  a $k$-graph $H$ on $n$ vertices \textit{$\varepsilon$-contains} $D^k_{n,s}$ if  there is a subset $U \subset V(H)$ with $|U|=k(s+1)-1$  such that $H$ $\varepsilon$-contains $D^k_{n,s}(U)$. A $k$-graph $H$ on $n$ vertices \textit{$\varepsilon$-contains} $H^k_{n,s}$ if $V(H)$ has a partition $V(H)=U \cup W$  with $|W|=s$ such that $H$ $\varepsilon$-contains $H^k_{n,s}(U,W)$.
Given $0<\theta<1$, we say a vertex $v\in V(H)$ is \emph{$\theta$-good} with respect to $H'$ if $|N_{H'} (v)\setminus N_H(v)|\leq \theta n^{k-1}$. Otherwise we say that $v$ is \emph{$\theta$-bad}. 
For a  $k$-graph $H$ and  $S\subseteq V(H)$, we use $H-S$ to denote the hypergraph obtained from $H$ by deleting $S$ and all edges of $H$ intersecting the set $S$, and we use $H[S]$ to denote the sub-hypergraph with vertex set $S$ and edge set $\{e\in E(H) : e\subseteq S\}$. For a  $k$-graph $H$ and  $E\subseteq E(H)$, we use $H-E$ to denote the hypergraph obtained from $H$ by removing edges from $E$. Let $K^k_n$ denote the complete $k$-graph on $n$ vertices.
By $x\ll y$ we mean that $x$ is sufficiently smaller than $y$ which needs to satisfy finitely many inequalities in the proof.
We omit the floor and ceiling functions when they do not affect the proof.

\section{Lower bounds on the anti-Ramsey number}
%
%
Recall inequality \eqref{lbhyper}. Note that if $G$ is an $s$-matching, then ${\cal G}=\{M_{s-1}\}$.
Therefore, we have the following lower bound for $\ar(n,k,M_s)$.
\begin{pro}\label{Lbound1}
	For $n\geq ks$, we have $\ar(n,k,M_s)\geq \ex(n,k,M_{s-1})+2$.
\end{pro}


Next we present a new construction for $n=ks$. Let $U$ be a subset of $V(K^k_n)$ such that $|U|=n-k-1$ and let $W=V(K^k_n)\setminus U$. Thus $|W|=k+1$. Let $f:E(K^k_n[U])\rightarrow [\binom{|U|}{k}]$ be a bijective coloring.

For an odd integer $k$, there are $\frac{1}{2}\binom{k+1}{(k+1)/2}$ distinct subsets $A_1,\ldots,A_{\frac{1}{2}\binom{k+1}{(k+1)/2}}$ of $W$ such that $|A_i|=(k+1)/2$ for $1\leq i\leq \frac{1}{2}\binom{k+1}{(k+1)/2}$ and $A_i\cap A_j\neq \emptyset$ for $1\leq i<j\leq \frac{1}{2}\binom{k+1}{(k+1)/2}$. Let $\mathcal{A}_i=\{e\in E(K^k_n):e\cap W=A_i \ \text{or} \  e\cap W= W\setminus A_i\}$ and  $\mathcal{H}_1$ be the complete $k$-graph $K^k_n$ with edge coloring $f_{\mathcal{H}_1}$, where
\begin{equation*}
	f_{\mathcal{H}_1}(e)=\left\{
	\begin{array}{ll}
		f(e), & \hbox{$e\in E(K^k_n[U])$;} \\
		\binom{|U|}{k}+i, & \hbox{$e\in \mathcal{A}_i$ for $1\leq i\leq \frac{1}{2}\binom{k+1}{(k+1)/2}$;} \\
		0, & \hbox{otherwise.}
	\end{array}
	\right.
\end{equation*}

For an even integer $k$, we fix a vertex $x\in W$. There are $\binom{k}{k/2-1}$ distinct subsets $B_1,\ldots,B_{\binom{k}{k/2-1}}$ of $W\setminus \{x\}$ such that $|B_i|=k/2-1$ for $1\leq i\leq \binom{k}{k/2-1}$. Let $\mathcal{B}_i=\{e\in E(K^k_n):x\in e \ \text{and} \ e\cap W=B_i\}\cup \{e\in E(K^k_n):e\cap W=W\setminus (B_i\cup \{x\})\}$ and let $\mathcal{H}_2$ be the $n$-vertex complete $k$-graph $K^k_n$ with edge coloring $f_{\mathcal{H}_2}$, where
\begin{equation*}
	f_{\mathcal{H}_2}(e)=\left\{
	\begin{array}{ll}
		f(e), & \hbox{$e\in E(\mathcal{H}_2[U])$;} \\
		\binom{|U|}{k}+i, & \hbox{$e\in \mathcal{B}_i$ for $1\leq i\leq \binom{k}{k/2-1}$;} \\
		0, & \hbox{otherwise.}
	\end{array}
	\right.
\end{equation*}
We can show the following lower bound on $\ar(n,k,M_{n/k})$ for $k \geq 3$. We remark that our construction here is new and different from the one for the graph case.
Note that for the case of $k=3$, the following proposition gives that $\ar(n,3,M_{n/3}) \geq \ex(n,3,M_{n/3-1})+5$, here $\ex(n,3,M_{n/3-1})=\binom{n-4}{3}$ by Theorem \ref{erdos3}.
\begin{pro}\label{Lbound2}
	$$\ar(n,k,M_{n/k})\geq\begin{cases}
		\binom{n-k-1}{k}+\frac{1}{2}\binom{k+1}{(k+1)/2}+2, & \text{$k$ is odd;}\\
		\binom{n-k-1}{k}+\binom{k}{k/2-1}+2, &\text{$k$ is even}.
	\end{cases}$$
\end{pro}
\begin{proof}
	It suffices to prove that neither $\mathcal{H}_1$ nor $\mathcal{H}_2$ has a rainbow perfect matching since the number of colors of $\mathcal{H}_1$ and $\mathcal{H}_2$ is $\binom{n-k-1}{k}+\frac{1}{2}\binom{k+1}{(k+1)/2}+1$ and $\binom{n-k-1}{k}+\binom{k}{k/2-1}+1$ respectively.
	
	For odd $k$, let $M$ be a matching in $\mathcal{H}_1$ covering $W$ such that $e\cap W\neq \emptyset$ for every $e\in M$. Note that $|M| \geq 2$ as $|W|=k+1$. If there is an edge $e_1\in M$ such that $|e_1\cap W|\neq (k+1)/2$, then there exists an edge $e_2\in M$ such that $e_2\neq e_1$ and $|e_2\cap W|\neq (k+1)/2$. Thus $f_{\mathcal{H}_1}(e_1)=f_{\mathcal{H}_1}(e_2)=0$ by the definition of $f_{\mathcal{H}_1}$ and $M$ is not a rainbow matching. If $|e\cap W|= (k+1)/2$ for  each edge $e\in M$, then the assumption that $M$ is a matching covering $W$ implies that
 $M=\{e_1,e_2\}$, where $e_1,e_2\in \mathcal{A}_i$ for some $i$.
The definition of $f_{\mathcal{H}_1}$ yields that $f_{\mathcal{H}_1}(e_1)=f_{\mathcal{H}_1}(e_2)$. Therefore, $M$ is not a rainbow matching.
	
	For even $k$,  let $M$ be a matching in $\mathcal{H}_2$ covering $W$ such that $e\cap W\neq \emptyset$ for every $e\in M$. Note that $|M| \geq 2$.  If there is no $e \in M$ such that $e  \in {\cal B}_i$ for some  $1 \leq i \leq \binom{k}{k/2-1}$, then   the definition of $f_{\mathcal{H}_2}$ implies that $f_{\mathcal{H}_2}(e_1)=f_{\mathcal{H}_2}(e_2)=0$ for any two distinct edges $e_1,e_2\in M$. Thus $M$ is not a rainbow matching. It remains to consider the case that there is an edge $g\in M$ such that $g\in \mathcal{B}_t$ for some $1\leq t\leq \binom{k}{k/2-1}$. As $ \mathcal{B}_i$ and $\mathcal{B}_t$ are cross-intersecting for each $1 \leq i \leq \binom{k}{k/2-1}$, then $M \cap \mathcal{B}_t=\{g\}$
and $M \cap \mathcal{B}_i=\emptyset$ for each $1 \leq i \neq t \leq \binom{k}{k/2-1}$, here we note that $ \mathcal{B}_i$ is intersecting and  $M$ is an matching.  As we assume $M$ covers $W$, there are at least
two edges $g_1, g_2 \in M \setminus\{g\}$ such that $g_1,g_2\notin \mathcal{B}_i$ for each $1 \leq i\leq \binom{k}{k/2-1}$. Thus 	$f_{\mathcal{H}_2}(g_1)=f_{\mathcal{H}_2}(g_2)=0$ and $M$ is not a rainbow matching.

\end{proof}

\section{A stability result in $3$-graphs}
In this section, we will prove the following stability result in 3-graphs which will be used in the proof of Theorem \ref{mainthm}.
\begin{lm}\label{mainsta}
	 Given reals $0<\varepsilon\ll c_0\ll 1$, there exists an integer $n_0$ such that the following holds. Let $H$ be a $3$-graph with $n>n_0$ vertices and $s$ be an integer. If $\nu(H)\leq s$ and
	\begin{equation}\label{cliedge}
		e(H)>\binom{3s+1}{3}+3s(n-3s-1),
	\end{equation}
	then the following holds.
	\begin{enumerate}[itemsep=0pt,parsep=0pt,label=$($\roman*$)$]
		\item For $5n/18-1\leq s\leq 13n/45$, if $H$ $\varepsilon$-contains $D^3_{n,s}$, then $H$ is a subgraph of $D^3_{n,s}$.
		\item For $13n/45\leq  s\leq (1-c_0)n/3$, $H$ is a subgraph of $D^3_{n,s}$.
	\end{enumerate}
\end{lm}

\noindent\textbf{Remark 1:}
Note that the condition (\ref{cliedge}) is tight as the following example. Define a hypergraph $\mathcal{D}(n,s)$ such that
\begin{equation*}
	E(\mathcal{D}(n,s))=\binom{[3s+1]}{3}\cup\left(\bigcup_{i=3s+2}^{n}\{\{1,i,x\}:2\leq x\leq 3s+1\}\right).
\end{equation*}
 One can observe that $\nu(\mathcal{D}(n,s))=s$, $e(\mathcal{D}(n,s))=\binom{3s+1}{3}+3s(n-3s-1)$ and $\mathcal{D}(n,s)$ is not a subgraph of $D^3_{n,s}$.

Before proving Lemma \ref{mainsta},  we recall the definition of the shifting.
Let $H$ be a $k$-graph on vertex set $[n]$. For vertices $1\leq i<j\leq n$, we define the \emph{$(i,j)$-shift} $S_{ij}$ by $S_{ij}(H)=\{S_{ij}(e):e\in E(H)\}$, where
	$$S_{ij}(e)=\begin{cases}
		e\setminus \{j\}\cup \{i\}, & \text{if $j\in e$, $i\notin e$ and $e\setminus \{j\}\cup \{i\}\notin E(H)$;}\\
		e, &\text{otherwise.}
	\end{cases}$$
	The following well-known result can be found in \cite{F95}.
	\begin{lm}\label{shiftprop}
For all $1\leq i<j\leq n$ and all $H$, the $(i,j)$-shift satisfies the following properties.
		\begin{enumerate}[itemsep=0pt,parsep=0pt,label=$($\roman*$)$]
			\item $e(H)=e(S_{ij}(H))$ and $|e|=|S_{ij}(e)|$,
			\item $\nu(S_{ij}(H))\leq \nu(H)$.
		\end{enumerate}
	\end{lm}
A $k$-graph $H$ is called \emph{stable} if $H=S_{ij}(H)$ for all $1\leq i<j\leq n$. It is not difficult to see that if $H$ is a stable $k$-graph, then for any subsets $\{u_1,\ldots,u_k\}, \{v_1,\ldots,v_k\}\subset[n]$ such that $u_i\leq v_i$ for $i\in[k]$, $\{v_1,\ldots,v_k\}\in E(H)$ implies $\{u_1,\ldots,u_k\}\in E(H)$. We need to introduce additional definitions.
Let $\omega(H)$ be the number of vertices in a largest complete subgraph of $H$. Two families $\mathcal{A}$ and  $\mathcal{B}$ are called \emph{cross-intersecting} if $A\cap B\neq\emptyset$ for all $A\in \mathcal{A}$, $B\in\mathcal{B}$.
Let $n,k,s$ be positive integers with $n\geq ks$ and $H$ be a $k$-graph on vertex set $[n]$. We say that $H$ is $s$-\emph{saturated}, if $\nu(H)\leq s$, but $\nu(\{e\}\cup H)=s+1$ for every $e\notin E(H)$.

Similar to Lemma 2 in \cite{LM}, we can prove the following one.
\begin{lm}\label{bigclique}
	There exist $\varepsilon>0$ and a positive integer $n_0$ such that the following holds. Let $H$ be an $s$-saturated $k$-graph on $n>n_0$ vertices and $s$ be an integer with $ n(1/k-1/2k^2)-1\leq s\leq (n-k+1)/k$. If  $(1-\varepsilon)ks\leq\omega(H)\leq ks+k-3$, then
	\begin{equation*}
		e(H)\leq\binom{ks+k-1}{k}-3.98\binom{(1-\varepsilon)ks}{k-1}+(2+8\varepsilon k^4)\binom{n}{k-1}.
	\end{equation*}
	
\end{lm}
As the proof of this lemma is an easy modification of the one for Lemma 2 in \cite{LM}, we include it in the appendix for the completeness. Relying on Lemma \ref{bigclique}, we can prove the following lemma.

\begin{lm}\label{cli-sta}
	There exist $\varepsilon>0$ and a positive integer $n_0$ such that the following holds. Let $H$ be a $3$-graph on $n>n_0$ vertices and $s$ be an integer with $ 5n/18-1\leq s\leq (n-2)/3$. If $\omega(H)\geq(1-\varepsilon)3s$, $\nu(H)\leq s$ and
	\begin{equation*}
		e(H)>\binom{3s+1}{3}+3s(n-3s-1),
	\end{equation*}
	then $H$ is a subgraph of $D^3_{n,s}$.
\end{lm}

\noindent\textbf{Remark 2:} Note that the condition for the number of edges is tight as the graph $\mathcal{D}(n,s)$ (see Remark 1.)

In order to prove Lemma \ref{cli-sta}, we need the following result.
\begin{lm}[Hilton and Milner \cite{HM}]\label{crossintersect}
	Suppose that $\mathcal{A},\mathcal{B}\subset \binom{[m]}{\ell}$ are non-empty and cross-intersecting, $m> 2\ell>0$. Then
	\begin{equation*}
		|\mathcal{A}|+|\mathcal{B}|\leq\binom{m}{\ell}-\binom{m-\ell}{\ell}+1.
	\end{equation*}
	Moreover, interchanging $\mathcal{A}$ and $\mathcal{B}$ and up to isomorphism, equality holds only in one of the following two cases:
	\begin{enumerate}[itemsep=0pt,parsep=0pt,label=$($\roman*$)$]
		\item $\ell=2$ and $\mathcal{A}=\mathcal{B}=\{\{1,x\}: 2\leq x\leq m\}$,
		\item $\mathcal{A}=\{[\ell]\}$,$\mathcal{B}=\{B\in\binom{[m]}{\ell}: B\cap[\ell]\neq\emptyset\}$.
	\end{enumerate}
\end{lm}

\noindent \textbf{Proof of Lemma \ref{cli-sta}.}
Note that we can turn a $3$-graph $H$ with $\nu(H)\leq s$ into an $s$-saturated $3$-graph $H'$ by adding edges.
 Thus we may assume that $H$ is  $s$-saturated and show $H$ is a subgraph of $D^3_{n,s}$.
	
	We claim that $\omega(H)\geq 3s+1$. Otherwise, if $\omega(H)\leq 3s$,  then Lemma \ref{bigclique} implies that
	\begin{equation*}
		e(H)\leq \binom{3s+2}{3}-3.98\binom{(1-\varepsilon)3s}{2}+(2+648\varepsilon)\binom{n}{2}.
	\end{equation*}
	Then, for $\varepsilon$ small enough and $n$  sufficiently large, we have
	\begin{align*}		
&e(H)-\binom{3s+1}{3}-3s(n-3s-1)\\ \leq&\binom{3s+2}{3}-3.98\binom{(1-\varepsilon)3s}{2}+(2+648\varepsilon)\binom{n}{2}-\binom{3s+1}{3}-3s(n-3s-1)\\ \leq&\binom{3s+2}{3}-3.9\binom{3s}{2}+2.1\binom{n}{2}-\binom{3s+2}{3}+\binom{3s+1}{2}-3s(n-3s-1)\\
			=&-3.9\binom{3s}{2}+2.1\binom{n}{2}+\binom{3s+1}{2}-3s(n-3s-1)\\
            =&-2.9\binom{3s}{2}+2.1\binom{n}{2}-3s(n-3s-2).
\end{align*}
 It is not difficult to check that $e(H)<\binom{3s+1}{3}+3s(n-3s-1)$ for $5n/18-1\leq s\leq (n-2)/3$ and $n$ large enough, a contradiction.
	
	For the case of $\omega(H)\geq 3s+2$, as $\nu(H)\leq s$, it follows that $H=K^3_{n,s}(U)$, where $U$ is a subset of $V(H)$ with $|U|=3s+2$.
	
	Now let us consider the case of $\omega(H)=3s+1$. Let $U$ be a largest clique in $H$ such that $|U|=3s+1$. Assume that $V(H)\setminus U=\{v_1,\ldots,v_{n-3s-1}\}$ and $F_i=N_H(v_i)$ for each $1 \leq i \leq n-3s-1$. We claim that $F_1,\ldots,F_{n-3s-1}\subseteq \binom{U}{2}$. Indeed, if there exists $F_i\nsubseteq \binom{U}{2}$, then there is an edge $e\in E(H)$ such that $v_i\in e$ and $|U\setminus e|\geq3s$. Thus $H[U\setminus e]$ contains a matching $M$ of size $s$ and then $M\cup \{e\}$ is a matching of size $s+1$ in $H$, a contradiction. Let $I=\{i: 1 \leq i \leq n-3s-1 \textrm{ and } F_i \not = \emptyset\}.$ If $|I|\leq 1$, then $H$ is a subgraph of $D^3_{n,s}$. Therefore, we may assume that $2\leq |I|\leq n-3s-1$.
We claim that $F_i$ and $F_j$ cross-intersecting for $i,j \in I$ and $F_i\neq F_j$. Otherwise, suppose that there are distinct $F_i$ and $F_j$ such that $F_i$ and $F_j$ are not cross-intersecting. Then there exist $f_i\in F_i$ and $f_j\in F_j$ such that $f_i\cap f_j=\emptyset$. Thus $M_1=\{f_i\cup\{v_i\},f_j\cup\{v_j\}\}$ is a matching of size two in $H$. Notice that $H[U\setminus V(M_1)]$ contains a matching $M_2$ of size $s-1$ as $|U\setminus V(M_1)|=3(s-1)$. Thus $M_1\cup M_2$ is a matching of size $s+1$ in $H$, a contradiction. By the same argument as above, if $F_i=F_j$ for some $i,j \in I$, then we get that $F_i$ is intersecting provided $|F_i| \geq 2$. In this case, we observe that $|F_i| \leq 3s$ by the famous Erd\H{o}s-Ko-Rado Theorem.
Let $I_1$ be a maximum  subset of $I$ such that $F_i \neq F_j$ for each $i,j \in I_1$.
By Lemma \ref{crossintersect}, we get that
	\begin{equation}
		\begin{split}
			\sum_{i \in I_1}|F_i|&=\frac{1}{|I_1|-1}\sum_{i,j \in I_1\atop i< j}(|F_i|+|F_j|)\\
			&\leq\frac{\binom{|I_1|}{2}} {|I_1|-1}\left(\binom{3s+1}{2}-\binom{3s-1}{2}+1\right)\\
			&\leq\frac{|I_1|}{2}\left(\binom{3s+1}{2}-\binom{3s-1}{2}+1\right)\\
			&=3s|I_1|.
		\end{split}
	\end{equation}
As we have shown $|F_i| \leq 3s$ for each $i \in I\setminus I_1$ and $|I| \leq n-3s-1$, we have
$$\sum_{i=1}^{n-3s-1}d_{H}(v_i)=\sum_{i \in I} |F_i| \leq 3s|I_1|+3s|I\setminus I_1| \leq 3s(n-3s-1).$$
 Therefore,
$e(H)=\binom{3s+1}{3}+\sum_{i=1}^{n-3s-1}d_{H}(v_i) \leq \binom{3s+1}{3}+3s(n-3s-1)$, which is a contradiction.
\qed

Let $S(H)$ be a graph obtained from $H$ by applying all possible $(i,j)$-shifts with $1 \leq i < j \leq n$.
Note that $S(H)$ is stable. We recall the following two results.


\begin{thm}[Gao, Lu, Ma and Yu \cite{GLMY}]\label{stablelemma}
	For any real $\varepsilon>0$, there exists a positive integer $n_1(\varepsilon)$ such that the following holds. Let $s,n$ be integers with $n\geq n_1(\varepsilon)$ and $1\leq s\leq n/3$, and let $H$ be a stable $3$-graph on the vertex set $[n]$. If $e(H)>\ex(n,3,M_s)-\varepsilon^4n^3$ and $\nu(H)<s$, then $H$ $\varepsilon$-contains $H^3_{n,s-1}([n]\setminus [s-1],[s-1])$ or $D^3_{n,s-1}([3(s-1)-1])$.
\end{thm}

\begin{thm}[Frankl \cite{F20}]\label{clique-sta}
Assume that $s \geq 27$.  Let $H$ be a $3$-graph on $3s$ vertices. If $\nu(H)< s$ and
	\begin{equation*}
		e(H)>\binom{3s-2}{3}+\binom{3s-2}{2}-\binom{3s-4}{2}+1,
	\end{equation*}
	then $H$ is a subgraph of $D^3_{n,s-1}$.
\end{thm}
 \noindent\textbf{Remark 3:} Note that Frankl \cite{F20} proved such a result for all $k \geq 3$. We only state a special case where $k=3$.


We are now ready to prove Lemma \ref{mainsta}.

\noindent
{\it Proof of Lemma \ref{mainsta}.}
To prove (i), notice that there exists a subset $U\subseteq V(H)$ of size $3s+2$ such that $|E(D^3_{n,s}(U))\setminus E(H)|\leq \varepsilon n^3$ since $H$ $\varepsilon$-contains $D^3_{n,s}$.
	Let $U=[3s+2]$ and $V(H)\setminus U=[n]\setminus [3s+2]$. Iterating the $(i,j)$-shift for all $1\leq i<j\leq n$ will eventually produce a  stable $3$-graph $S(H)$. By Lemma \ref{shiftprop}, we have $e(S(H))=e(H)$ and $\nu(S(H))\leq \nu(H)\leq s$. By the definition of $(i,j)$-shift, $|E(D^3_{n,s}(U))\setminus E(S(H))|\leq \varepsilon n^3$. We claim that there is a complete subgraph of size at least $3(1-3\varepsilon^{1/3})s$ in $S(H)$. Let $U'=[3(1-3\varepsilon^{1/3})s]$. Suppose that $S(H)[U']$ is not a complete subgraph, then $U\setminus U'$ is an independent set in $S(H)$ as $S(H)$ is stable. Thus $|E(D^3_{n,s}(U))\setminus E(S(H))|\geq \binom{|U\setminus U'|}{3}\geq\binom{5\varepsilon^{1/3}n/2}{3}> \varepsilon n^3$ for sufficiently large $n$, a contradiction.
 Lemma \ref{cli-sta} implies that $S(H)$ is a subgraph of $D^3_{n,s}$.
	Suppose that $H$ is not a subgraph of $D^3_{n,s}$. Then in the process of producing $S(H)$ by $(i,j)$-shift, there
is a $3$-graph $H'$ obtained from $H$ by the series of shifts such that $H'$ has exactly $n-3s-3$ isolated vertices. By Lemma \ref{shiftprop}, it follows that $\nu(H')\leq s$ and $e(H')=e(H)$. Note that
   \begin{equation*}
			e(H)>\binom{3s+1}{3}+3s(n-3s-1)
			\geq \binom{3s+1}{3}+\binom{3s+1}{2}-\binom{3s-1}{2}+1.
	\end{equation*}
 Theorem \ref{clique-sta} gives that $H'$ is a subgraph of $D^3_{n,s}$, a contradiction. Thus $H$ is a subgraph of $D^3_{n,s}$.

To prove (ii), through the $(i,j)$-shift for all $1\leq i<j\leq n$, we obtain a stable $3$-graph $S(H)$ such that $e(S(H))=e(H)$ and $\nu(S(H))\leq \nu(H)\leq s$. For  $13n/45\leq  s\leq (1-c_0)n/3$, note that $e(H)>\binom{3s+1}{3}+3s(n-3s-1)\geq \ex(n,3,s+1)-\varepsilon^4n^3$.
 Theorem \ref{stablelemma} implies that $S(H)$ either $\varepsilon$-contains $D^3_{n,s}([3s-1])$ or $\varepsilon$-contains $H^3_{n,s}([n]\setminus [s],[s])$. Let $W=[s]$ and $U=[n] \setminus W$. We next show that the latter case is impossible.
Suppose that $S(H)$ $\varepsilon$-contains $H^3_{n,s}(U,W)$.
Let $T=\{v\in W : d_{S(H)}(v)<\binom{n-1}{2}-\sqrt{\varepsilon} n^2\}$. We claim $|T|\leq 3\sqrt{\varepsilon}n$. Otherwise,
\begin{equation}\label{1}
	|E(H^3_{n,s}(U,W))\setminus E(S(H))|> 3\sqrt{\varepsilon}n\cdot\sqrt{\varepsilon} n^2/3=\varepsilon n^3,
\end{equation}
a contradiction.
 If we set $s=\alpha n$, then $13/45\leq \alpha \leq (1-c_0)/3$ as $13n/45\leq  s\leq (1-c_0)n/3$. Let $f(x)=\frac{9x^3}{2}-\frac{1-(1-x)^3}{6}$. We get that
\begin{equation}\label{diff2}
	\begin{split}
	e(S(H)[U])&>\binom{3s+1}{3}+3s(n-3s-1)-\left(\binom{n}{3}-\binom{n-s}{3}\right)\\
	&=\frac{(1-(1-\alpha)^3)n^3}{6}-\frac{9\alpha^3n^3}{2}+o(n^{3})=f(\alpha)n^3+o(n^3).
	\end{split}
\end{equation}
Since $f'(x)=\frac{26x^2+2x-1}{2}$ is increasing in $[13/45,(1-c_0)/3]$ with $f'(13/45)>0$, we have $f(\alpha)\geq f(13/45)>0.001$ for $13/45\leq \alpha \leq (1-c_0)/3$. By inequality (\ref{diff2}), we have $e(S(H)[U])>10\sqrt{\varepsilon}n^3$ for sufficiently small $\varepsilon$. Let $U'=\{s+1,\ldots,s+10\sqrt{\varepsilon}n\}$. We claim that $S(H)[U']$ is a clique. Otherwise, as $S(H)$ is stable, we get that $U \setminus U'$ is an independent set. It follows that $e\cap  U' \neq \emptyset$ for each $e\in E(S(H)[U])$. We can see that $e(S(H)[U])\leq 10\sqrt{\varepsilon}n^3$, a contradiction to the lower bound on $e(S(H)[U])$.
Let $M$ be a matching of size $|U'|/3$ in $S(H)[U']$ and  $T'=W \setminus T$.  Thus $|T'|\geq s-3\sqrt{\varepsilon}n$ as $|T|\leq 3\sqrt{\varepsilon}n$. Let $M'$ be a maximum matching in $S(H)-V(M)$ such that $|e\cap T'|=1$ for each $e\in M'$. We claim that $|M'|=|T'|$. Otherwise, suppose that $|M'|<|T'|$. Then there exists a vertex $v\in T'\setminus V(M')$ such that $N_{S(H)}(v)\cap \binom{V(H)\setminus (V(M)\cup V(M')\cup T')}{2}=\emptyset$. Thus $d_{S(H)}(v)\leq \binom{n-1}{2}-\binom{n-|V(M)\cup V(M')\cup T'|}{2}\leq \binom{n-1}{2}-\binom{(c_0-10\sqrt{\varepsilon})n}{2}<\binom{n-1}{2}-\sqrt{\varepsilon} n^2$, a contradiction.
Notice that $|M|=|U'|/3>3 \sqrt{\varepsilon} n$ and then $M\cup M'$ is a matching of size at least $s+1$, a contradiction to the condition $\nu(S(H))\leq \nu(H)\leq s$. Thus $S(H)$ $\varepsilon$-contains $D^3_{n,s}$.
Similar to the proof of (i), there is a complete subgraph of size at least $3(1-3\varepsilon^{1/3})s$ in $S(H)$. Thus $S(H)$ is a subgraph of $D^3_{n,s}$ by Lemma \ref{cli-sta}. Repeating arguments in the proof of (i), we can show that $H$ is a subgraph of $D^3_{n,s}$. \hfill $\square$

\section{Proof of Theorem \ref{mainthm}}

To prove Theorem \ref{mainthm}, we need the following result.
\begin{thm}[Guo, Lu and Mao \cite{GLM}]\label{stability}
	Let $\varepsilon,\rho$ be two reals such that $0<\rho\ll\varepsilon<1$. Let $n,s$ be two integers such that $n$ is sufficiently large and $n/54+1\leq s\leq 13n/45+1$. Let $H$ be a $3$-graph on vertex set $[n]$. If $e(H)>\ex(n,3,M_s)-\rho n^3$ and $\nu(H)\leq s-1$, then $H$  $\varepsilon$-contains $H^3_{n,s-1}$ or $D^3_{n,s-1}$.
\end{thm}

We first prove the following lemma.
\begin{lm}\label{per-sta}
	Let $\varepsilon,\beta$ be constants such that $0<\varepsilon< 10^{-6}$ and $3\sqrt{\varepsilon}<\beta<1$. Let $n$ be a sufficiently large integer and $H$ be a $3$-graph on vertex set $[n]$. If $e(H)>(1-\varepsilon)\binom{n}{3}$ and $d_H(v)>\beta n^2$ for every $v\in[n]$, then $H$ has a matching covering all but at most two vertices.
\end{lm}

\begin{proof}
 The assumption that $e(H)>(1-\varepsilon)\binom{n}{3}$ gives that all but at most $\sqrt{\varepsilon }n$ vertices in $H$ have degree at least $(1-\sqrt{\varepsilon })\binom{n-1}{2}$. Otherwise,
	\begin{equation*}
			\binom{n}{3}-e(H)>\frac{1}{3}\left(\sqrt{\varepsilon } n\cdot \sqrt{\varepsilon } \binom{n-1}{2}\right)=\varepsilon \binom{n}{3},
	\end{equation*}
	a contradiction.
	
	Let $R=\{v\in V(H):d_H(v)<(1-\sqrt{\varepsilon })\binom{n-1}{2}\}$ and $r=|R|$. Denote the vertices in $R$ by $v_1,\ldots,v_{r}$ where $r\leq \sqrt{\varepsilon }n$. One can greedily find a matching $\{e_1,e_2,...,e_r\}$ in $H$ such that $e_i\cap R=\{v_i\}$ for all $1\leq i\leq r$. Suppose $\{e_1,e_2,...,e_t\}$ is a matching in $H$ with $t<r$ satisfying the condition above.  Since $d_{H}(v)>\beta n^2>3rn$ for $v\in V(H)$, there exists an edge $e_{t+1}\in E(H)$ such that $e_{t+1}\cap R=\{v_{t+1}\}$ and $e_{t+1}\cap(\cup_{i=1}^t e_i))=\emptyset$. Continuing the process, we may find the desired matching $M=\{e_1,\ldots,e_{r}\}$. 	
	
	Let $H'=H-V(M)$ and $n'=|V(H')|$. For every vertex $x\in V(H')$, it satisfies that $d_{H'}(x)>(1-\sqrt{\varepsilon })\binom{n}{2}-|V(M_1)|n>(1-7\sqrt{\varepsilon} )\binom{n'}{2}>\binom{n'-1}{2}-\binom{\lceil2n'/3\rceil}{2}$ for sufficiently large $n$. By Theorem \ref{degree-version}, $H'$ has a matching $M'$ covering all but at most two vertices in $V(H')$. Thus $M\cup M'$ is a matching covering all but at most two vertices in $V(H)$.
\end{proof}

We prove the following upper bound on $\ar(n,3,M_s)$ for relatively small $s$.
\begin{lm}\label{upper1}
	 For a given real $0<c_0\ll 1$, there exists an integer $n_0=n_0(c_0)$ such that $\ar(n,3,M_s)\leq \ex(n,3,M_{s-1})+2$ for $n/6\leq s\leq (1-c_0)n/3$ and $n>n_0$.
\end{lm}

\begin{proof}
Note that $\ex(n,3,M_{s-1})=\max\{\binom{n}{3}-\binom{n-s+2}{3},\binom{3s-4}{3}\}$ by  Theorem \ref{erdos3}. Let $c(n,s)=\max\{\binom{n}{3}-\binom{n-s+2}{3},\binom{3s-4}{3}\}+2$
and $f_{n,s}:E(K^3_n)\rightarrow [c(n,s)]$ be a surjective coloring. We use $H$ to denote the edge-colored $K^3_n$.

  Define $G$ as a subgraph of $H$ with $c(n,s)$ edges such that each color appears on exactly one edge of $G$. Let $\varepsilon, \rho$ be reals such that $0<\rho\ll\varepsilon\ll c_0$. Note that $e(G)=\max\{\binom{n}{3}-\binom{n-s+2}{3},\binom{3s-4}{3}\}+2\geq \max\{\binom{n}{3}-\binom{n-s+1}{3},\binom{3s-1}{3}\}-\rho n^3$ for sufficiently large $n$.
  Therefore, if $n/6\leq s\leq 13n/45+1$, then $G$ either $\varepsilon$-contains $H^3_{n,s-1}$ or $\varepsilon$-contains  $D^3_{n,s-1}$ by Theorem \ref{stability}. The proof is split into the following three cases.

\medskip
	\textbf{Case 1.}  $n/6\leq s\leq 13n/45+1$ and  $G$  $\varepsilon$-contains $H^3_{n,s-1}$.
\medskip

	Since $G$ $\varepsilon$-contains $H^3_{n,s-1}$, there is a partition $U,W$ of $V(G)$ such that $|U|=n-s+1$, $|W|=s-1$ and
	\begin{equation}\label{1}
		|E(H^3_{n,s-1}(U,W))\backslash E(G)|\leq\varepsilon n^3.
	\end{equation}
	Let $T'=\{v\in W\ |\ d_G(v)<\binom{n-1}{2}-\sqrt{\varepsilon} n^2\}$.
Fix $x\in W$ and define
\begin{align*}
 T=\left\{
   \begin{array}{ll}
     T', & \hbox{if $T'\neq \emptyset$;} \\
     \{x\}, & \hbox{otherwsie.}
   \end{array}
 \right.
\end{align*}
Let $t=|T|$. We claim $1\leq t\leq 3\sqrt{\varepsilon}n$, otherwise,
\begin{equation}\label{1}
		|E(H^3_{n,s-1}(U,W))\backslash E(G)|> 3\sqrt{\varepsilon}n\cdot\sqrt{\varepsilon} n^2/3=\varepsilon n^3,
	\end{equation}
a contradiction.

Let $W'=W\setminus T$ and $G'=G-W'$.	
One can see that
\begin{align*}
e(G')&\geq\binom{n}{3}-\binom{n-s+2}{3}+2-\left({n\choose 3}-{n-|W'|\choose 3}\right)\\
&={n-|W'|\choose 3}-\binom{(n-|W'|)-(t-1)}{3}+2.
\end{align*}
By Theorem \ref{2sk}, $G'$ contains a rainbow matching of size $t+1$. Let $M$ be a rainbow $(t+1)$-matching in $G'$. Since $d_G(v)\geq\binom{n-1}{2}-\sqrt{\varepsilon} n^2$ for every $v\in W'$, we can greedily construct a matching $M'$ of size $s-t-1$ in
$G-V(M)$ such that $|e\cap W'|=1$ and $f_{n,s}(M)\cap f_{n,s}(M')=\emptyset$. More precisely, let $W'=\{v_1,\ldots,v_{s-t-1}\}$.
	For $v_1\in W'$, note that
	$d_G(v_1)>\binom{n-1}{2}-\sqrt{\varepsilon} n^2\geq \binom{n-1}{2}-\binom{n-|V(M)\cup W'|}{2}+s$.
Since $G$ is a rainbow graph, there exists an edge $e_1\in E(G-V(M))$ such that $e_1\cap W'=\{v_1\}$ and $f_{n,s}(e_1)\notin f_{n,s}(M)$.
	Now suppose we have found a rainbow matching
	$\{e_1,e_2,...,e_r\}$ in $G-V(M_1)$  such that
	$|e_i\cap W'|=\{v_i\}$ and $f_{n,s}(e_i)\notin f_{n,s}(M_1)$ for all $i\in [r]$.
	If $r=s-t-1$, then $M\cup \{e_1,\ldots, e_{s-t-1}\}$ is a desired matching.
	So we may assume that $r<s-t-1$.
	Write $G_r=G-V(M)-(\cup_{i=1}^r e_i)$. Note that $|[n]\setminus (W'\cup V(M)\cup(\cup_{i=1}^r e_i))|\geq2n/15$.
	Since $d_G(v_{r+1})>\binom{n-1}{2}-\sqrt{\varepsilon} n^2\geq \binom{n-1}{2}-\binom{n-|W'\cup V(M)\cup(\cup_{i=1}^r e_i)|}{2}+s$,  there exists an edge $e_{r+1}\in E(G_r)$ such that
	$e_{r+1}\cap W'=\{v_{r+1}\}$ and $f_{n,s}(e)\notin f_{n,s}(M_1)\cup\{f_{n,s}(e_1),\ldots,f_{n,s}(e_r)\}$. Continuing the process, we are able to find the desired rainbow matching $M'=\{e_1,\ldots,e_{s-t-1}\}$.
	Now $M\cup M'$ is a rainbow $s$-matching in $H$.

	%

%


\medskip
\textbf{Case 2.} $n/6\leq s\leq 5n/18$ and $G$ $\varepsilon$-contains $D^3_{n,s-1}$.

\medskip

In this case, we write $s=\alpha n$ with $1/6\leq \alpha \leq 5/18$. Let $f(x)=\frac{1-(1-x)^3}{6}-\frac{9x^3}{2}$. Then the following holds:
\begin{equation}\label{diff}
		\binom{n}{3}-\binom{n-s+2}{3}-\binom{3s-4}{3}=\frac{(1-(1-\alpha)^3)n^3}{6}-\frac{9\alpha^3n^3}{2}+o(n^{3})=f(\alpha)n^3+o(n^3).
\end{equation}
Since $f'(x)=\frac{1-2x-26x^2}{2}$ is decreasing in $[1/6,5/18]$ with $f'(1/6)<0$, we have $f(\alpha)\geq f(5/18)>0.007$ for $1/6\leq \alpha \leq 5/18$.

As $G$ $\varepsilon$-contains $D^3_{n,s-1}$, there is a subset $U$ of size $3s-4$ such that $|E(D^3_{n,s-1}(U))\setminus E(G)|\leq \varepsilon n^3$. We claim that there is a subset $S\subseteq  V(G)\setminus U$ such that $|S|=\varepsilon^{1/6} n/4$ and $d_{G}(v)\geq \varepsilon^{1/6} n^{2}$ for each $v\in S$. Otherwise, if there are at most $\varepsilon^{1/6}n/4$ vertices with degree at least $\varepsilon^{1/6}n^2$, then the number of edges intersecting $V(G)\setminus U$ is no more than $\varepsilon^{1/6}n^3/4+\varepsilon^{1/6}n^3= 5\varepsilon^{1/6}n^3/4.$
Thus by inequality (\ref{diff}), we have
$$e(G)\leq \binom{3s-4}{3}+5\varepsilon^{1/6}n^3/4<\binom{n}{3}-\binom{n-s+2}{3},$$ a contradiction.

Let $M$ be a maximum matching in $G$ such that $|e\cap S|=1$ for every $e\in M$. We claim that $|M|=|S|$. Otherwise, $|M|<|S|= \varepsilon^{1/6} n/4$ implies that there exists a vertex $v\in S\setminus V(M)$ such that $N_{G}(v)\cap \binom{V(G)\setminus (V(M)\cup S)}{2}=\emptyset$. Thus $d_G(v)\leq |V(M_1)\cup S|\cdot n\leq 3\varepsilon^{1/6} n^{2}/4$, a contradiction to  the fact that $d_G(v)\geq \varepsilon^{1/6} n^{2}$.

Recall that $|E(D^3_{n,s-1}(U))\setminus E(G)|\leq \varepsilon n^3$. Then all but at most $3\sqrt{\varepsilon}n$ vertices in $G$ are $\sqrt{\varepsilon}$-good. Otherwise,
\begin{equation*}
		|E(D^3_{n,s-1}(U))\setminus E(G)|=\frac{1}{3}\sum_{v\in V(G)}|N_{D^3_{n,s-1}(U)}(v)\backslash N_{G}(v)|
		>(3\sqrt{\varepsilon}n\cdot \sqrt{\varepsilon}n^{2})/3=\varepsilon n^3,
\end{equation*}
a contradiction. Let $U^{bad}$ be the set of all $\sqrt{\varepsilon}$-bad vertices. Note that $|U^{bad}|\leq 3\sqrt{\varepsilon}n$. Let $U'=U\setminus (U^{bad}\cup V(M))$.
 Recall $|M|<\varepsilon^{1/6} n/4$. Then we get $|U^{bad}\cup V(M)|\leq3\sqrt{\varepsilon}n+3\varepsilon^{1/6} n/4\leq 2\varepsilon^{1/6}n$ and $|U'|>|U|-3\sqrt{\varepsilon}n-2\varepsilon^{1/6} n/4\geq |U|-3\varepsilon^{1/6}n/4$.
Therefore,
$$e(G[U'])\geq e(D^3_{n,s-1}(U))-|E(D^3_{n,s-1}(U))\setminus E(G)|-|U^{bad}\cup V(M_1)|\cdot\binom{U}{2}\geq (1-3\varepsilon^{1/6})\binom{|U'|}{3}.$$
In addition,
$$d_{G[U']}(x)\geq \binom{|U|}{2}-\sqrt{\varepsilon}n^2-|U^{bad}\cup V(M_1)|\cdot |U|\geq 1/2|U'|^2$$ for each $x\in U'$.
Notice that $s-|M|<|U'|/3$. By Lemma \ref{per-sta}, $G[U']$ has a matching $M'$ of size $s-|M|$ and then $M\cup M'$ is a rainbow $s$-matching.

\medskip
\textbf{Case 3.} $5n/18+1\leq s\leq 13n/45+1$ and $G$  $\varepsilon$-contains $D^3_{n,s-1}$ or $13n/45+2\leq s\leq (1-c_0)n/3$.

\medskip

 Theorem \ref{erdos3} guarantees the existence of a rainbow $(s-1)$-matching in $G$.  Let $\{M_1,\ldots,M_t\}$ be a maximum family of edge-disjoint rainbow $(s-1)$-matchings in $G$ such that $f_{n,s}(M_i)\cap f_{n,s}(M_j)=\emptyset$ for all $i\neq j$.
Let $S_i=[n]\setminus V(M_i)$ for all $1\leq i\leq t$ and let $\ell=|S_1|=\cdots=|S_t|$. The assumption $s\leq (1-c_0)n/3$ implies that $\ell\geq c_0n+3$. We next show $t \leq s$. We need only to consider the case of $t\geq 2$. In this case, if $|S_i\cap S_j|\geq 3$ for $i\neq j$, then there exists an edge $e\subseteq S_i\cap S_j$ such that $f_{n,s}(e)\notin f_{n,s}(M_i)$ or $f_{n,s}(e) \notin f_{n,s}(M_j)$. Therefore, either $\{e\}\cup M_i$ or $\{e\}\cup M_j$ is a rainbow $s$-matching. It remains to consider the case where $|S_i\cap S_j|\leq 2$ for all $i\neq j$. Thus $\binom{S_i}{3}\cap\binom{S_j}{3}=\emptyset$ for all $i\neq j$. Since $\cup_{i=1}^t\binom{S_i}{3}\subseteq \binom{{n}}{3}$, we have $|\cup_{i=1}^t\binom{S_i}{3}|=\sum_{i=1}^t|\binom{S_i}{3}|=t\binom{\ell}{3}\leq\binom{n}{3}$. As $\ell\geq c_0n+3$, we get that $t\leq\binom{n}{3}/\binom{\ell}{3}\leq n/4<s$. Let $G_1$ be a hypergraph obtained from $G$ by removing all edges in $M_1\cup\cdots\cup M_t$.

For sufficiently large $n$, it holds that
\begin{equation*}
	e(G_1)=e(G)-t(s-1)
	\geq \binom{3s-4}{3}+2-s(s-1)
	>\binom{3s-5}{3}+2+(3s-6)(n-3s+5).
\end{equation*}
 Note that $G_1$ is a subgraph of $D^3_{n,s-2}(U)$ for some  $U$ of size $3s-4$ by Lemma \ref{mainsta}.


Recall $H$ is the edge-colored $K_n^3$. Let $H_1=H-E(H[U])$.
Next we show that $H_1$ has a rainbow $2$-matching $M$ such that $|V(M)\cap U|\leq 2$.
Let $e_1\in E(H_1)$ be an edge contained in $[n]\setminus U$ and $P=\{e\in E(H_1):e\cap e_1=\emptyset\}$. If there exists an edge $g_1\in P$ such that $f_{n,s}(e_1)\neq f_{n,s}(g_1)$, then $\{e_1,g_1\}$ is a rainbow $2$-matching such that $|(e_1\cup g_1)\cap U|\leq 2$. Thus we assume that edges in $P$ receive the same color as $e_1$. Since the number of colors used for coloring $E(H)$ is at least $\binom{3s-4}{3}+2$ and $e(H_1)=e(H)-\binom{3s-4}{3}$, we get that $H_1$ has at least two edges with distinct colors. Thus there is an edge $g_2\in E(H_1)$ such that $e_1\cap g_2\neq \emptyset$ and $f_{n,s}(e_1)\neq f_{n,s}(g_2)$. Since $n\geq c_0n+3s$, there exists an edge $e_3$ such that $e_3\cap (e_1\cup g_2)=\emptyset$ and $e_3\subseteq[n]\setminus U $. Note that $e_3 \in P$ and  $f_{n,s}(e_3)=f_{n,s}(e_1)$. Therefore,
$\{e_3,g_2\}$ is a rainbow $2$-matching such that $|(e_3\cup g_2)\cap U|\leq 2$.

Let $M=\{e,g\}$ be a rainbow $2$-matching in $H_1$ such that $|V(M)\cap U|\leq 2$. Note that $|U \setminus V(M)|\geq 3s-6$. We choose $U_1$ as a subset of $U \setminus V(M)$ with $|U_1|=3s-6$. Let $G_2$ be a subgraph of $G_1[U_1]$ by removing  edges colored by $f_{n,s}(e)$ and $f_{n,s}(g)$. Notice that we only remove two edges from $G_1[U_1]$. Recall that $G_1$ is a subgraph of $D_{n,s-2}^3(U)$ and edges in $G_1$ receive distinct colors. Additionally, the number edges in $G_1$ which contains a vertex from $U\setminus U_1$ is at most $\binom{3s-5}{2}+\binom{3s-6}{2}$. Since $e(G_1)\geq \binom{3s-4}{3}+2-s(s-1)$, we have
\begin{align*}
	e(G_2) &\geq \binom{3s-4}{3}+2-s(s-1)-\binom{3s-5}{2}-\binom{3s-6}{2}-2\\
           &=\binom{3s-6}{3}-s(s-1) \\
           &>\binom{3s-6}{3}+1-\binom{3s-7}{2} \\
           & =\binom{3s-7}{3}+1,
\end{align*}
here we note that $s$ is large as we assume $n$ is large and $s \geq 13n/45+2$.
Thus by Theorem \ref{erdos3}, there is a rainbow $(s-2)$-matching $M'$ in $G_2$ such that $f_{n,s}(M)\cap f_{n,s}(M')=\emptyset$. Notice that $M\cup M'$ is a rainbow $s$-matching in $H$ as desired.
\end{proof}

\begin{lm}\label{upper2}
	 For a given real $0<c_0\ll 1$, there exists an integer $n_0=n_0(c_0)$ such that for $n>n_0$,
\begin{equation*}
 \ar(n,3,M_s)\leq  \left\{
    \begin{array}{ll}
      \ex(n,3,M_{s-1})+2, & \hbox{ if $(1-c_0)n/3\leq s< n/3$ ;} \\
      \ex(n,3,M_{s-1})+5, & \hbox{ if $s=n/3$.}
    \end{array}
  \right.
\end{equation*}
\end{lm}

	\begin{proof}
Notice that $\ex(n,3,M_{s-1})=\binom{3s-4}{3}$ for $(1-c_0)n/3\leq s \leq n/3$.
Let \begin{equation*}
 c(n,s)= \left\{
    \begin{array}{ll}
      {3s-4\choose 3}+2, & \hbox{ if $(1-c_0)n/3\leq s< n/3$ ;} \\
      {3s-4\choose 3}+5, & \hbox{ if $s=n/3$.}
    \end{array}
  \right.
\end{equation*}
For a surjective mapping $f_{n,s}:E(K^3_n)\rightarrow [c(n,s)]$, where $V(K^3_n)=[n]$, we use  $H$ to denote the edge-colored $K^3_n$. Let $G$ be a subgraph of $H$ with $c(n,s)$ edges such that each color appears on exactly one edge of $G$. We assume that $d_G(1)\geq d_G(2)\geq \cdots \geq d_G(n)$ without losing any generality. Let $U=[3s-4]$ and $W=[n]\setminus U$. We define $R=\{x\in U:d_{G[U]}(x)<n^2/15\}$ and $r=|R|$. Then
			\begin{equation}\label{eq:37-1}
			\binom{3s-4}{3}-e(G[U])>\frac{1}{3}r\left(\binom{3s-5}{2}-\frac{n^2}{15}\right).
		\end{equation}
		Let $H'=H-E(H[U\setminus R])$.

\medskip
		\textbf{Claim 1.} For $n$ large enough, $r<2c_0n$.
\medskip

\noindent
{\it Proof of Claim 1:}
Recall that $n$ is large and  $(1-c_0)n/3\leq s \leq n/3$.
For each $x\in R$, we have
 \[
d_G(x)<n^2/15+n(n-(3s-4))\leq n^2/15+n(c_0n+4)<n^2/14.
 \]
Thus $d_G(y)<n^2/14$ for every $y\in[n]\setminus U$. The fact  $|[n]\setminus U|\leq c_0n+4$ implies that
\begin{align}\label{eq:37-2}
e(G[U])> \binom{3s-4}{3}+2-|[n]\setminus U| n^2/14\geq\binom{3s-4}{3}+2-c_0n^3/7.
\end{align}
Combining (\ref{eq:37-1}) and  (\ref{eq:37-2}), we can get that $r<2c_0n$. This completes the proof of Claim 1.  	

\medskip
		\textbf{Claim 2.} If $H'$ has a rainbow matching $M$ such that $|V(M)\cap (W\cup R)|\geq r+4$, then $G$ has a rainbow matching of size $s$.
\medskip

\noindent
{\it Proof of Claim 2:}
Let $M$ be a rainbow matching of $H'$ such that $|V(M)\cap (W\cup R)|\geq r+4$. We may assume that $|M|\leq r+4$. It follows that
\begin{align}\label{eq:37-3}
 |V(M)\cap (U\setminus R)|\leq 3|M|-4-r\leq 2r+8.
\end{align}

Then
\[
|U\setminus (V(M)\cup R)|=3s-4-r-|V(M)\cap (U\setminus R)|\geq 3(s-|M|).
 \]
Let $P=\{e\in E(G):f_{n,s}(e)\in f_{n,s}(M)\}$ and  $G'=G[U\setminus(V(M)\cup R)]-P$. We show that there exists a matching of size $s-|M|$ in $G'$. 
Recall that $r<2c_0n$. Then
\begin{equation*}
	\begin{split}
	e(G')>&e(G[U])-|(V(M)\cup R)\cap U|n^2-|M|\\
	>&\binom{3s-4}{3}+2-\frac{c_0n^3}{7}-(3r+12)n^2-(r+4)\quad \mbox{ (\ref{eq:37-2}) and (\ref{eq:37-3})}\\
	>&(1-\beta)\binom{|V(G')|}{3},
	\end{split}
\end{equation*}
where $0<\beta\ll 1$. Recall that  $d_{G[U]}(x)>n^2/15$ for every vertex $x\in U\setminus R$. Thus for every $x\in V(G')$, it satisfies that
\begin{align*}
  d_{G'}(x)&\geq d_{G[U]}(x)-(|U|-1)|(V(M_1)\cap U)\cup R|-(r+4)\\
&\geq n^2/15-(3s-5)(3r+8)-(r+4) >|V(G')|^2/20.
\end{align*}
 By Lemma \ref{per-sta}, $G'$ has a matching $M'$ of size $s-|M|$  since $|V(G')|=|U\setminus (V(M)\cup R)|\geq 3(s-|M|)$. Then $M\cup M'$ is a rainbow $s$-matching of $H$.  This completes the proof of Claim 2.

By Claim 2, it suffices to show that $H'$ has a rainbow matching $M$ such that $|V(M)\cap (W\cup R)|\geq r+4$.
There are two cases.
		
		\medskip
		\textbf{Case 1.} $r>0$.
     	\medskip

		Let $G_1=G-E(G[U\setminus R])$. Since $e(G)\geq\binom{3s-4}{3}+2$ and $e(G[U\setminus R])\leq \binom{3s-4-r}{3}$, we have
\begin{align}\label{eq:37-4}
e(G_1)\geq \binom{3s-4}{3}-\binom{3s-4-r}{3}\geq r\binom{3s-4-r}{2}.
\end{align}
 We claim that there are at least $r+4$ vertices $x_1,\ldots,x_{r+4}\in (V(G_1)\setminus U)\cup R$ such that $d_{G_1}(x_i)>3(r+4)n$. Otherwise we have
	    \begin{equation*}
			\begin{split}
			e(G_1)&\leq 3(r+4)n\bigg(n-(3s-4-r)-(r+3)\bigg)+\frac{n^2}{14}(r+3)\\
				&\leq3(r+4)c_0n^2+3(r+4)n+\frac{n^2}{14}(r+3)\quad \mbox{(Claim 1 and  $c_0\ll 1$)}\\
				&\leq r\binom{(1-4c_0)n}{2}\\
				&<r\binom{3s-4-r}{2},
			\end{split}
		\end{equation*}
 a contradiction to $(\ref{eq:37-4})$.
	
	Now we greedily find a matching $M$ such that every edge in $M$ contains exactly one vertex from $\{x_1,\ldots,x_{r+4}\}$. Let $S=\{x_1,\ldots,x_{r+4}\}$. Since $d_{G_1}(x_1)>3(r+4)n$, there exists one edge $e_1$ in $G_1$ such that $e_1\cap S=\{x_1\}$.  Now suppose that we have found a matching $\{e_1,e_2,...,e_t\}$ in $G_1$ such that $e_i\cap S=\{x_i\}$ for all $1\leq i\leq t$. For $t<r+4$, note that  $d_{G_1}(x_i)>3(r+4)n$ and $|S\cup (\cup_{i=1}^te_i)|\leq r+4+2t$. Thus there exists an edge $e_{t+1}\in E(G_1)$ such that $e_{t+1}\cap S=\{v_{t+1}\}$ and $e_{t+1}\cap(\cup_{i=1}^t e_i)=\emptyset$. Continuing the process, we manage to find a desired matching $M=\{e_1,\ldots,e_{r+4}\}$.

	\medskip
\textbf{Case 2.} $r=0$.
  \medskip

We distinguish the following two subcases.

\medskip
\textbf{Subcase 2.1.} $n>3s$.
\medskip

It suffices to show that $H'$ has a rainbow $2$-matching $M$ such that $|V(M)\cap W|\geq 4$. Note that $|W|= |[n]-U|= n-(3s-4)\geq 5$.
We begin to consider the case that every edge in ${W\choose 3}$ is colored by the same color.
Recall that $H'=H-E(H[U\setminus R])$.
Thus the number of colors used for edges in $E(H')$ is at least $\binom{3s-4}{2}+2-\binom{|U|}{3}=2$. So there exists an edge $e\in E(H')$ such that  $1\leq |e\cap W|\leq 2$ and $f_{n,s}(e)\notin f_{n,s}({W\choose 3})$. One can see  $|W\setminus e|\geq 3$. Hence we may choose   $e'\in {W\setminus e\choose 3}$. Then $\{e,e'\}$ is a desired rainbow matching.

We next assume that there exist two edges $e_1,e_2\in {W\choose 3}$ such that $f_{n,s}(e_1)\neq f_{n,s}(e_2)$. If $e_1\cap e_2=\emptyset$, then $\{e_1,e_2\}$ is a desired rainbow matching. If $|e_1\cap e_2|=2$, since $|W|\geq 5$,  then we are able to pick an edge $e_3\in E(H')$ such that $e_3\cap (e_1\cup e_2)=\emptyset$. Then either $\{e_1,e_3\}$ or $\{e_1,e_2\}$ is a desired rainbow matching. If $|e_1\cap e_2|=1$, then there exists $g\in E(H')$ such that $|g\cap e_2|=2$ and $|g\cap e_1|=2$. One can see that $f_{n,s}(e_1)\neq f_{n,s}(g)$ or $f_{n,s}(e_2)\neq f_{n,s}(g)$. Then by the previous argument, we may find a desired rainbow 2-matching.

\medskip
\textbf{Subcase 2.2.} $n=3s$.
\medskip

In this subcase, notice that $|W|=4$ and $|f_{n,s}(E(H')) |\geq 5$. It suffices to prove that there exists a rainbow matching covering $W$.

\medskip
\textbf{Claim 3.} Let $S\in {W\choose 1}\cup {W\choose 2}$. If there exist two distinct edges $e_1,e_2\in E(H')$ such that $S= e_1\cap W= e_2\cap W$, and $f_{n,s}(e_1)\neq f_{n,s}(e_2)$, then we are able to find a rainbow matching $M$ covering $W$.
\medskip

\noindent
{\it Proof of Claim 3:}
Let $e_1$ and $e_2$ be two edges satisfying the  condition and $S'=W\setminus S$. Let $g\in E(H')$ such that $g\cap W=S'$ and $g\cap (e_1\cup e_2)=\emptyset$. Since $f_{n,s}(e_1)\neq f_{n,s}(e_2)$, we get that either $\{e_1,g\}$ or $\{e_2,g\}$ is a rainbow matching covering $W$.
This completes the proof of Claim 3.

 Write $W=\{x_1,x_2,x_3,x_4\}$.  For $ S\in {W\choose 1}\cup {W\choose 2}$, let $B_S=\{f_{n,s}(e)\ |\ S=e\cap W, e\in E(H')\}$.
Note that $|f_{n,s}(E(H')) |\geq 5$. We have the following Claim.

\medskip
\textbf{Claim 4.}  The following holds:
 \begin{itemize}
   \item [(i)] $|B_S|=1$ for any  $ S\in {W\choose 1}\cup {W\choose 2}$;

   \item [(ii)] $B_S=B_{W\backslash S}$ for any $S\in {W\choose 1}\cup{W\choose 2}$;

   \item [(iii)] $|\cup_{S\in {W\choose 1}\cup {W\choose 2} }B_S|\geq 5$.
 \end{itemize}

 \noindent
 {\it Proof of Claim 4:} Claim 3 implies that (i) is true. For (ii), if there is an $S \in  {W\choose 1}\cup{W\choose 2}$ such that  $B_S \not =B_{W\backslash S}$, then we are able to find two disjoint edges $e_1,e_2 \in H'$ such that $S=e_1 \cap W$, $W \setminus S=e_2 \cap W$, and $f_{n,s}(e_1) \neq f_{n,s}(e_2)$. Thus  lemma \ref{upper2} follows from Claim 2 and we are done.  The fact that $|f_{n,s}(E(H')) |\geq 5$ tells us that (iii) is true.

For simplicity, by (i), we may write  $\{b_{i,j}\}=B_S$ if $S=\{x_i,x_j\}$ and $\{b_{i}\}=B_S$ if $S=\{x_i\}$.

\medskip
\textbf{Claim 5.} There exists a 2-set $\{x_i,x_j\}\in {W\choose 2}$  such that $|\{b_{i,j},b_{i},b_{j}\}|=3$.

\noindent
{\it Proof of Claim 5:} One can observe that $|\{b_{i,j} \ |\ \{i,j\}\in {[4]\choose 2}\}|\leq 3$ by  (ii). Thus by (iii), we have $|\{b_{i}\ |\ i\in [4]\}|\geq 2$ and there exists $\{x_i,x_j\}\in {W\choose 2}$ such that $b_{i,j}\notin \{b_{\ell}\ |\ \ell\in [4]\}$.
Assume that $b_{1,2}\notin \{b_{\ell}\ |\ \ell\in [4]\}$. If either $b_{1}\neq b_{2}$ or $b_{3}\neq b_{4}$, then
either $\{x_1,x_2\}$ or $\{x_3,x_4\}$ is a desired 2-set. We assume that $b_{1}= b_{2}$ and $b_{3}= b_{4}$.
Then  by (ii) and (iii), $|\{b_{\ell}\ |\ \ell\in [4]\}\cup \{b_{1,i}\ |\ i\in \{2,3,4\}\}|\geq 5$, we have $|\{b_{1,2},b_{1,3},b_{1,4}, b_{1},b_{3}\}|=5$. Then
$\{b_{1,3}, b_{1},b_{3}\}$ is a desired 2-set. This completes the proof of Claim 5.

By Claim 5, without loss of generality,   we may assume that $|\{b_{1,2},b_{1},b_{2}\}|=3$. Recall that $b_{1,2}=b_{3,4}$. Let  $e_1,e_2,e_3$ be three pairwise disjoint edges such that $\{x_3,x_4\}=e_1\cap W$, $\{x_1\}=e_2\cap W$, and $\{x_2\}=e_3\cap W$. Then $\{e_1,e_2,e_3\}$ is a rainbow matching covering $W$. This completes the proof of the lemma.
\end{proof}
We are ready to present the proof of Theorem \ref{mainthm}.

\noindent \textbf{Proof of Theorem \ref{mainthm}.}
By Theorem \ref{2sk}, we need only to show $\ar(n,3,M_s)= \ex(n,3,M_{s-1})+2$ for $n/6\leq s < n/3$ and $\ar(n,3,M_s)= \ex(n,3,M_{s-1})+5$ for $n=3s$. For $n/6\leq s < n/3$, the lower bound follows from Proposition \ref{Lbound1} while the upper bound follows from Lemmas \ref{upper1} and \ref{upper2}.
  For $n=3s$, the combination of Proposition \ref{Lbound2} and
  Lemma \ref{upper2} implies that $\ar(n,3,M_s)= \ex(n,3,M_{s-1})+5$.

\qed

\section{Concluding remarks}
Let $k\geq 3$, $s\geq 3$ and $n=ks$.  \"Ozkahya and Young \cite{OY} conjectured that
$$\ar(n,k,M_s)=\begin{cases}
		\ex(n,k,M_{s-1})+2, & \text{if $s<c_k$};\\
		\ex(n,k,M_{s-1})+k+1, &\text{if $s\geq c_k$},
	\end{cases}$$
where $c_k$ is a constant depending on $k$.

By Proposition \ref{Lbound2}, we have
$$\ar(n,k,M_s)\geq\begin{cases}
	\binom{ks-k-1}{k}+\frac{1}{2}\binom{k+1}{(k+1)/2}+2, & \text{$k$ is odd;}\\
	\binom{ks-k-1}{k}+\binom{k}{k/2-1}+2, &\text{$k$ is even}.
\end{cases}$$

If $n=ks$, then $\ex(n,k,M_{s-1})=\binom{ks-k-1}{k}$ by a result from \cite{Fr172}. Note that $\tfrac{1}{2}\binom{k+1}{(k+1)/2}+2 >k+1$ for $k \geq 3$.
Thus Conjecture \ref{mainconj} is not true for $n=ks$. Although we showed that our construction indeed gives the true value of $\ar(n,3,M_{n/3})$ for $n$ large enough,  we do not know whether it is the case for $k \geq 4$.

{\bf Appendix.}

\noindent
{\bf Proof of Lemma \ref{bigclique}.}
Let $H=(V,E)$ and $U$ be a largest clique of $H$. Note that $(1-\varepsilon)ks\leq |U|\leq ks+k-3$.
We choose $M$ as a matching of size $s$ in $H$ which maximizes $|V(M)\cup U|$ and set $M'=\{e\in M: e\not\subseteq U\}$. For $n$ large enough, the following claim was proved explicitly in the proof of Lemma 2 in \cite{LM}.

\medskip
\noindent
\textbf{Claim A.} \\
(i) $|V(M)\cup U|=ks+k-1$. \\
			(ii) $|M'|\leq 2\varepsilon ks$.\\
			(iii) Each edge of $H$ either is contained in $U$, or intersects an edge of $M'$.

Let $H'$ denote the $k$-graph with vertex set $V(H')=V$ and edge set $E(H')=\binom{V(M)\cup U}{k}$. Clearly, the size of the largest matching in $H'$ is $s$. We say that a subset $f\subset V$ of $\ell$ vertices is \emph{thick} if it is contained in more than $3\varepsilon k^3\binom{|U|}{k-\ell}$ edges $e\in E$ with $e\subseteq U\cup f$, and \emph{thin} otherwise. The following claim follows from the proof of Lemma 2 in \cite{LM} explicitly.

\medskip
\noindent
\textbf{Claim B.}
		If a subset $f$ of $\ell$ elements is thick, then each $k$-element subset of $U\cup f$ containing $f$ is an edge of $H$.

	Now let $A=V(M)\setminus U$ and $a=|A|$. Observe first that every vertex in $A$ is thin. Indeed, if $\{w\}\subseteq A$ is thick, then all $k$-element subsets of $U\cup \{w\}$ belong to $H$ by Claim 6. Thus $U\cup\{w\}$ is a larger clique than $U$, which is a contradiction to the choice of $U$.  Using this fact, one can show the following:
	\begin{equation}\label{E'-E}
		|E(H')\setminus E(H)|\geq (1-3\varepsilon k^3)a\binom{|U|}{k-1}.
	\end{equation}
	
	Now we estimate the number of edges in $E(H)\setminus E(H')$. Let $E_1$ be the set of edges which have at least two vertices in $V(M')$. We define $E_2$ to be the set of all edges $e$ such that $e\cap V(M')=\{w\}$ and the set $(e\setminus U)\cup \{w\}$ is thin. Similarly, let $E_3$ be the set of all edges $e$ such that $e\cap V(M')=\{w\}$ and the set $(e\setminus U)\cup \{w\}$ is thick. By the proof of Lemma 2 in \cite{LM}, we have
	\begin{equation}\label{E1}
		|E_1|\leq \varepsilon k^4a\binom{n}{k-1},
	\end{equation}

\begin{equation}\label{E2}	
		|E_2|\leq 3\varepsilon k^4 a(\binom{n}{k-1}-\binom{|U|}{k-1}),	
\end{equation}

	\begin{equation}\label{E3}
			|E_3|\leq a(\binom{n}{k-1}-0.99\binom{|U|}{k-1}).
	\end{equation}
	Therefore, from inequalities (\ref{E'-E}), (\ref{E1}), (\ref{E2}), and (\ref{E3}), we get
	\begin{equation}\label{more2}
			e(H')-e(H)\geq a(1.99\binom{|U|}{k-1}-\binom{n}{k-1}-4\varepsilon k^4\binom{n}{k-1}).
	\end{equation}

	Note that $a\geq 2$ since $\omega(H)\leq ks+k-3$. Thus
	\begin{equation*}
		\begin{split}		
			e(H)\leq  &\binom{ks+k-1}{k}-a(1.99\binom{|U|}{k-1}-\binom{n}{k-1}-4\varepsilon k^4\binom{n}{k-1})\\
			\leq &\binom{ks+k-1}{k}-3.98\binom{(1-\varepsilon)ks}{k-1}+(2+8\varepsilon k^4)\binom{n}{k-1}.
		\end{split}	
	\end{equation*}

\end{document}